\documentclass[11pt,reqno]{amsart}

\usepackage{amsfonts, amssymb, amsthm, amsmath, bbm, wasysym, enumerate, amsxtra, latexsym, fullpage, amscd, graphics, epic, sansmath, mathtools, multicol, url, hyperref, bbm, blindtext, wasysym, comment, datetime, stmaryrd}
\usepackage[shortlabels]{enumitem}
\usepackage[all,cmtip]{xy}
\usepackage[dvipsnames]{xcolor}
\usepackage[makeroom]{cancel}  
\usepackage[enableskew]{youngtab}
\usepackage{tikz}
\usetikzlibrary{cd}

\usepgflibrary{shapes.geometric}

\tikzstyle{V}=[draw, fill =black, circle, inner sep=0pt, minimum size=1.5pt]
\tikzstyle{C}=[draw, fill =white, circle, inner sep=0pt, minimum size=1.5pt]
\tikzstyle{over}=[draw=white,double=black,line width=2pt, double distance=.5pt]
\usetikzlibrary{arrows, positioning, calc, chains, cd}
\tikzset{
	ch/.style={circle,draw,on chain,inner sep=2pt},
	chj/.style={ch,join},
	every path/.style={shorten >=4pt,shorten <=4pt}
	}
\newcommand{\dnode}[2][chj]{%
	\node[#1,label={below:#2}] (#1) {};}
\newcommand{\dnodenj}[1]{%
	\dnode[ch]{#1}}
\newcommand{\dydots}{%
	\node[chj,draw=none,inner sep=1pt] {\dots};}
\numberwithin{equation}{section}
\usepackage[margin = 0.9in]{geometry}
\theoremstyle{definition}

\newtheorem{theorem}{Theorem}[section]
\newtheorem{lemma}[theorem]{Lemma}

\newtheorem{ex}[theorem]{Example}

\newtheorem{rmk}[theorem]{Remark}

\def\<{\langle}
\def\>{\rangle}

\newcommand{\A}{\mathrm{A}}
\newcommand{\B}{\mathrm{B}}

\newcommand{\CC}{\mathbb{C}}

\newcommand{\cD}{\mathcal{D}}
\newcommand{\cH}{\mathcal{H}}

\newcommand{\inv}{^{-1}}

\newcommand{\ld}{\lambda}

\newcommand{\rStd}{\textup{rStd}}

\newcommand{\Sh}{\textup{Sh}}

\newcommand{\Std}{\textup{Std}}

\newcommand{\tif}{\textup{if }}

\newcommand{\uM}{\underline{M}}

\newcommand{\ZZ}{\mathbb{Z}}

\allowdisplaybreaks 

\newcommand{\Bc}[1]{     \begin{cases}   #1   \end{cases}}
\newenvironment{eq}{\begin{equation}}{\end{equation}}

\title{A short proof on the transition matrix from the Specht basis to the Kazhdan--Lusztig basis}
\author[Mee Seong Im]{Mee Seong Im}
\address{Department of Mathematical Sciences, United States Military Academy, West Point, NY 10996}
    \email{meeseongim@gmail.com}

 
\begin{document}

\begin{abstract}
We provide a short proof on the change-of-basis coefficients from the Specht basis to the Kazhdan--Lusztig basis, 
using Kazhdan--Lusztig theory for parabolic Hecke algebra. 
\end{abstract}
 
\maketitle 

\section{Introduction}
It is well-known that the irreducible representations for the symmetric group $\Sigma_{n}$ are the Specht modules parametrized by the set of partitions of $n$, or equivalently, set of standard Young tableaux consisting of $n$ boxes. 
The Specht module has a purely combinatorial basis (called Specht basis) described in terms of certain alternating sums of so-called tabloids.

For the Young tableaux consisting of two rows of equal size, Russell--Tymoczko in \cite{RT17} compare the Specht basis with another combinatorial basis (called the web basis) which arises from Temperley--Lieb algebra and knot theory. 
In their context, the web basis is a reincarnation of the Kazhdan--Lusztig basis, but this is not true in general.
Their main result is a combinatorial model that gives a different proof to a special case of a classic theorem by Naruse \cite[Theorem~4.1]{Na89} that the change-of-basis matrix is unitriangular, altogether with some vanishing conditions. The unitriangularity result is also found in \cite{GM88}, but without vanishing conditions.

In this paper, we give a short proof of \cite[Theorem~4.1]{Na89} using Kazhdan--Lusztig theory for parabolic Hecke algebra.
We also notice that the argument applies to all classical types, with a nonstandard notion of tableaux. 
For example, in type B we use certain centro-symmetric tableaux which are not the bi-tableaux that parametrize the irreducibles. 
That is to say, the analog of the Specht modules here are not the irreducibles.
However, we hope that they correspond to the top cohomology of the Springer fiber corresponding to the parabolic subgroup.
Furthermore, we omit the details for type D as future work. 
In this manuscript, we introduce a Specht basis for this module using certain alternating sum of centro-symmetric tabloids. 
We then prove a unitriangular theorem regarding the change-of-basis matrix between the Specht basis and the Kazhdan--Lusztig basis (cf. Theorem~\ref{thm:main1}).

\section{Combinatorics}
\subsection{Weyl groups}
\begin{figure}[ht!]
\label{figure:Dynkin}
\caption{Dynkin diagrams of type $\A_{d}$ and $\B_d$.}
\label{figure:Dyn}
\[
\begin{array}{ccc}
\begin{tikzpicture}[start chain]
\dnode{$1$}
\dnode{$2$}
\dydots
\dnode{$d$}
\end{tikzpicture}
&
\quad
\begin{tikzpicture}[start chain]
\dnode{$0$}
\dnodenj{$1$}
\dydots
\dnode{$d-1$}
\path (chain-1) -- node[anchor=mid] {\(=\joinrel=\joinrel=\)} (chain-2);
\end{tikzpicture}
\\
\textup{Type }\A_d \quad (d \geq 1)
&
\textup{Type }\B_d \quad (d \geq 2)
\end{array}
\]
\end{figure}
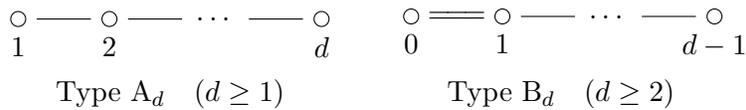

For type $\Phi = \A$ or $\B$, let $W^\Phi_d$ be the Weyl group of type $\Phi_d$ associated with the Dynkin diagrams in Figure~\ref{figure:Dynkin}. 

Write 
\eq
[d] := [1, d] \cap \ZZ,
\quad
[\pm d] := [-d, d] \cap \ZZ.
\endeq
Set
\eq
I^\A_d = [d] 
\quad 
\mbox{ and }
\quad 
I^\B_d 
=  
[0, d-1] \cap \ZZ.
\endeq

We denote by $S^\Phi_d = \{s^\Phi_i : i \in I^\Phi_d\}$, the corresponding generators of $W^\Phi_d$ as Coxeter groups.
Denote by $\ell^\Phi:W^\Phi_d \to \ZZ_{\geq 0}$ the length function.

For any set $I$, 
let $\textup{Perm}(I)$ be the group of permutations on $I$, and let $(i,j) \in \textup{Perm}(I)$ be the transposition for $i,j \in I$.
It is a standard fact that $W^\Phi_d$ can be identified as a group of certain permutations (see \cite{BB05}). 
Here we treat them as subgroups of $\textup{Perm}([\pm d])$ as follows:
\eq
\begin{split}
W^\B_d &:=  \{g \in \textup{Perm}([\pm d]) :  g(-i) = -g(i) \textup{ for all } i\},
\\
W^\A_{d-1} &:=  \{g \in W^\B_d : \textup{neg}(g) = 0\},
\end{split}
\endeq
where the function $\textup{neg}$ counts the total number of negative entries of $\{1, \ldots, d\}$, i.e.,
\eq\label{def:neg}
\textup{neg}: W^\B_d \to \ZZ_{\geq 0},
\quad g \mapsto \#(\{ k \in [d] : g(k) < 0\}).
\endeq
The generators can be identified with the following permutations in $[\pm d]$:
\eq
\begin{split}
s_i &= s_i^\A = s_i^\B =  (i, i+1)(-i,-i-1) \quad\textup{for}\quad 1 \leq i \leq d-1,
\\
s_0^\B &= (-1, 1). 
\end{split}
\endeq

It is sometimes convenient to use the one-line notation: 
\eq\label{def:perm}
w  \equiv |w(1), w(2), \ldots , w(d)| 
\quad 
\mbox{ for } 
w \in W^\Phi_d.
\endeq
We remark that the element $w$ is uniquely determined by these $d$ values due to the centro-symmetry condition:  $g(-i) = -g(i)$ for all $i$.

\subsection{Parabolic {H}ecke algebras}
Let $\cH^\Phi = \cH(W^\Phi_d)$ be the Hecke algebra of $W^\Phi_d$ over $\CC[q,q\inv]$. As a free module, $\cH^\Phi$ has a basis $\{ H^\Phi_w : w \in W^\Phi_d\}$.
The multiplication in $\cH^{\Phi}$ is determined by
\eq\label{def:Hecke}
\begin{split}
H^\Phi_w H^\Phi_x &= H^\Phi_{wx} \quad \tif \ell^\Phi(wx) = \ell^\Phi(w) + \ell^\Phi(x),
\\
(H^\Phi_s)^2 &= (q\inv-q) H^\Phi_s + H^\Phi_e \quad \tif s\in S^\Phi,   
\end{split}
\endeq
where $e\in W^\Phi_d$ is the identity, and hence $H^\Phi_e$ is the identity element in $\cH^\Phi$.
By a slight abuse of notation, we use the same symbol $H_i$ $(1\leq i \leq n-1)$ to denote the element $H^\Phi_{s^\Phi_i} \in \cH^\Phi$. 
We write $H^\Phi_0 = H^\Phi_{s^\Phi_0}$ for $\Phi = \B$.
It is a standard fact that $\cH^\Phi$ is generated as a $\CC[q,q\inv]$-algebra by $\{H^\Phi_i : i \in I^\Phi_d\}$.

For each subset $J \subset I^\Phi_d$, denote the corresponding  parabolic subgroup and parabolic Hecke algebra of $W^\Phi_d$ by
\eq
W^\Phi_J = \<s_j^\Phi : j\in J\> 
\quad \mbox{ and } \quad 
\cH(W^\Phi_J) = \<H_j^\Phi : j\in J\>, 
\endeq
respectively. 
We also denote the set of shortest right coset representatives for $W^\Phi_J\setminus W^\Phi_d$ by
\eq\label{def:DJ}
\cD^\Phi_J 
= \{ w\in W^\Phi_d : \ell^\Phi(wg) = \ell^\Phi(w) + \ell^\Phi(g) \textup{~for all~} w\in W^\Phi_J\}.
\endeq
Denote the induced trivial module corresponding to $J \subset I^\Phi_d$ by
\eq
M^J_\Phi = \CC[q,q\inv] \otimes_{\cH_J^\Phi} \cH^{\Phi},
\endeq
where $\CC[q,q\inv]$ is regarded as a right $\cH^{\Phi}$-module by setting
\eq
f \cdot H^\Phi_i = q\inv f 
\quad
\textup{for}
\quad
f \in \CC[q,q\inv].
\endeq
The module $M^J_\Phi$ admits a standard basis
\eq
\{M_w = 1\otimes H^\Phi_w : w \in \cD^\Phi_J\}.
\endeq
We note that $M^J_\Phi$ can also be identified as the right $\cH^{\Phi}$-module $y_\ld\cH^{\Phi}$ under the assignment $M_w \mapsto y_J H^\Phi_w$, where $y_J = \sum_{w\in W^\Phi_J} q^{-\ell(w)} H^\Phi_w$, on which $\cH^{\Phi}$ acts by a right multiplication.
The action of $\cH^{\Phi}$ on $M_\Phi^J$ below follows from \eqref{def:Hecke} and a straight-forward calculation:
\[
M_w \cdot H^\Phi_i =
\begin{cases}
M_{w s_i} &\tif w s_i \in \cD^\Phi_J, \quad \ell^\Phi(ws_i) > \ell^\Phi(w), 
\\ 
M_{w s_i} +(q\inv - q) M_w &\tif w s_i \in \cD^\Phi_J, \quad  \ell^\Phi(ws_i) < \ell^\Phi(w), 
\\
q\inv M_w &\tif w s_i \not\in \cD^\Phi_J.
\end{cases}
\]
Following \cite[Theorem~3.1]{Soe97}, $M_\Phi^J$ admits a Kazhdan--Lusztig basis $\{\uM_w : w \in \cD^\Phi_J\}$ such that
\eq
\uM_w = \sum_{x \in \cD^\Phi_J} m^\Phi_{x,w}(q\inv) M_x,
\endeq
where $m^\Phi_{x,w} \in \ZZ[q\inv]$.
We also define polynomials $p^\Phi_{w,x}(q\inv)\in \ZZ[q\inv]$ such that
\eq
M_w = \sum_{x\in \cD^\Phi_J} p^\Phi_{w,x}(q\inv) \uM_x.
\endeq
Below we recall some properties of the (inverse) parabolic Kazhdan--Lusztig polynomials.
\lemma\label{lem:PKL}
We have the following: 
\begin{enumerate}[(a)]
\item if $x \not\leq w$ with respect to the Bruhat order on $W^\Phi_d$, then $m^\Phi_{x,w} = 0 = p^\Phi_{x,w}$, 
\item $m^\Phi_{x,x} = 1 = p^\Phi_{x,x}$.
\end{enumerate}
\endlemma
\proof
It follows from \cite[Theorem~3.1]{Soe97} that the matrix $(m_{x,w}(q\inv))_{x,w}$ is upper unitriangular, its inverse matrix, $(p_{x,w}(q\inv))_{x,w}$, is also upper unitriangular. 
\endproof
\subsection{Young tableaux}
Now we want to associate each parabolic subgroup $W^\Phi_J$ of $W^\Phi_d$ a ``Young tableau'' beyond type A.
The idea here is to use the corresponding composition 
(cf. \cite{DJ86, BKLW18}). 
To be precise, we set
\eq
\Lambda^\A(d) =\bigsqcup_{n \geq 0} \Lambda^\A(n,d) 
\quad 
\mbox{ and }
\quad 
\Lambda^\B(d) =\bigsqcup_{n \geq 0} \Lambda^\B(n,d),
\endeq
where the legit $n$-part (signed) compositions are defined as 
\begin{align}
\Lambda^\A(n,d) &=
\left\{
\lambda = (\lambda_i)_{i\in [n]} \in \ZZ_{>0}^{n}
: 
\sum_{i=1}^n \lambda_i = d
\right\},
\\
\Lambda^\B(n,d) &= 
\Bc{
\left\{
(\lambda_i)_{i\in [\pm r]} \in \ZZ_{>0}^{n} : 
\begin{array}{c}
\lambda_0 \in 1+2\ZZ, 
\:\:
\displaystyle{\sum_{i=-r}^r} \ld_i= 2d+1,
\\
\lambda_{-i} = \lambda_i
~\forall i
\end{array}\right\} 
&\tif n = 2r+1, 
\\
\left\{
\ld \in \Lambda^\B(2r+1,d) : 
\ld_0 = 1
\right\} 
&\tif n = 2r. 
}
\end{align}
In other words, the corresponding parabolic subgroups are generated by  
\eq
S^\A(d) - \left\{ s_{\lambda_1}, s_{\lambda_1+\lambda_2}, \ldots, s_{\lambda_1+\ldots +s_{\lambda_{n-1}}}\right\} 
\quad \mbox{ for } \quad  \ld \in \Lambda^\A(n,d), 
\endeq 
and 
generated by (abbreviate $\lfloor {\frac{\lambda_0}{2}}\rfloor$ by $\lambda^\natural_0$ here)
\eq
S^\B(d) - 
\begin{cases}\left\{s_{\lambda_1}, s_{\lambda_1+\lambda_2}, \ldots, s_{\lambda_1 + \ldots +\lambda_{r-1}}\right\}
&\tif n=2r, 
\\
\left\{s_{\lambda^\natural_0}, s_{\lambda^\natural_0+\lambda_1}, \ldots, s_{\lambda^\natural_0 +\lambda_1 + \ldots +\lambda_{r-1}}\right\}
&\tif n= 2r+1, 
\end{cases}
\quad \mbox{ for } \quad  \lambda \in \Lambda^\B(n,d). 
\endeq
We also write $W^\Phi_\ld$ and $\cH^\Phi_\ld$  
to denote the corresponding parabolic subgroups and the parabolic Hecke algebras, respectively.

For type A, any composition $\ld = (\ld_1, \ldots, \ld_n) \in \Lambda^\A(n,d)$ defines a Young subgroup
\[
W^\A_\ld \simeq \Sigma_{\ld_1} \times \Sigma_{\ld_2} \times \ldots \times \Sigma_{\ld_n}.
\]
That is, any two compositions giving the same Young subgroup (and hence parabolic Hecke algebra) correspond to the same partition.
We denote the set of type A partitions by
$
\Pi^\A(d) = \bigsqcup_{n \geq 1} \Pi^\A(n,d),
$
where
\eq
 \Pi^\A(n,d) = \{ \ld \in \Lambda^\A(n,d) : \ld_1 \geq \ld_2 \geq \ldots \geq \ld_n\}.
\endeq
However, for type B, a composition $\ld = (\ld_0, \ldots, \ld_n) \in \Lambda^\B(n,d)$ defines a Young subgroup
\[
W^\B_\ld \simeq W^\B_{\ld^\natural_0} \times \Sigma_{\ld_1}\times \Sigma_{\ld_2} \times \ldots \times \Sigma_{\ld_n}.
\]
In particular, $W^\B_\ld$ is always a product of symmetric groups when $n$ is even.
Hence, compositions $\ld,\mu$ describe the same Young subgroup (and hence parabolic Hecke algebra) if $\ld_0 = \mu_0$ and also that $(\ld_1, \ldots, \ld_n)$ corresponds to the same partition as $(\mu_1, \ldots, \mu_n)$.
We denote the set of type B ``partitions'' by
$
\Pi^\B(d) = \bigsqcup_{n \geq 1} \Pi^\B(n,d),
$
where
\eq
 \Pi^\B(n,d) = \{ \ld \in \Lambda^\B(n,d) : \ld_1 \geq \ld_2 \geq \ldots \geq \ld_n\}.
\endeq  

\ex
Let $d=3$. There are 8 subsets of $I^\B_d$, and hence we have
\[
\begin{array}{c|ccccccccc}
\ld \in \Lambda^\B(3)& W^\B_J & J
\\
\hline
(7) &W^\B_3 & \{0,1,2\}
\\
(1,5,1) &\<s_0^\B, s_1\> \simeq W^\B_2 & \{0,1\}
\\
(2,3,2)&\<s_0^\B, s_2\> \simeq \Sigma_2 \times \Sigma_2 & \{0,2\}
\\
(3,1,3)&\<s_1, s_2\> \simeq \Sigma_3  & \{1,2\}
\\
(1,1,3,1,1)&\<s_0^\B\> \simeq \Sigma_2 & \{0\}
\\
(1,2,1,2,1)&\<s_1\> \simeq \Sigma_2 &\{1\}
\\
(2,1,1,1,2)&\<s_2\> \simeq \Sigma_2 &\{2\}
\\
(1,1,1,1,1,1,1)& \{e\} & \varnothing 
\end{array}
\]
\endex

Let $\ld \in \Lambda^\Phi(n,d)$.
For now we assume that $\Phi = \A$. 
Let its corresponding Young diagram be
$[\ld] = \{(i,j) : i\geq 1, 1\leq j \leq \ld_i\}$.
A Young tableau of shape $\ld$ (or $\ld$-tableau for short) is a bijection
\[
T: [\ld] \to [d].
\]
For example, when $\ld = (m,m)$ or $(m,m,m)$, a $\ld$-tableau represents a $2\times m$ or $3\times m$ grid, respectively, whose boxes are filled in by every number, exactly once, from 1 to $d=2m$ or $d=3m$, respectively. 
We also write $\Sh(T) = \ld$.
The set of all young tableaux admits a (right) action of the symmetric group $W^\A_{d-1}=\Sigma_d$ by permuting the letters that are filled in the boxes. 

A Young tableau is called row standard if the entries in each row are increasing; it is called standard if the entries in each row and each column are increasing. Denote by $\rStd(\ld)$ and $\Std(\ld)$ the sets of row standard and standard Young tableaux of shape $\ld$, respectively. 
We remark that the size of $\Std(\ld)$ is counted by the hook formula, and hence
\eq
\#\Std((m,m)) = \frac{1}{m+1}{2m\choose m} 
\quad 
\mbox{ and }
\quad 
\#\Std((m,m,m)) = \frac{2 \cdot (3m)!}{m! \cdot (m+1)! \cdot (m+2)!}. 
\endeq
It is well-known that the following assignment is a bijection:
\eq
\rStd(\ld) \to \cD^\A_\ld,
\quad T \mapsto w_T,
\endeq 
where $w_T$ is given by row-reading of tableau $T$.

We now generalize this bijection to type B by defining the Young tableaux of classical type. For $\Phi = \B$,
 $\ld \in \Lambda^\Phi(n,d)$, denote the corresponding Young diagram by
\eq
[\ld] = \begin{cases}
\left\{(0,j) : -\ld_1 \leq j \leq \ld_1\right\}
\sqcup \left\{\pm(i,j) : i\geq 1, 0\leq j \leq \ld_i-1\right\}
&\tif n = 2r+1, 
\\
\left\{\pm(i,j) : i\geq 1, 1\leq j \leq \ld_i \right\} &\tif n=2r.
\end{cases}
\endeq
Note that when $d$ is even, the corresponding Young diagram is just a type A Young diagram, with a copy obtained by rotating it 180 degrees. 
\ex
For $d=2$, $\Pi^\B(2) = \{(2,1,2), (1,1,1,1,1)\}$, and hence the tableaux are
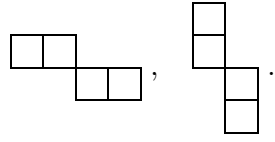
\[
\young(~~::,::~~)~,
\quad
\young(~:,~:,:~,:~)~.
\]
For $d=3$, we then have
\[
\begin{array}{c|ccccccccc}
\ld \in \Pi^\B(3)
&(7) 
&(1,5,1)
&(2,3,2)&
(3,1,3)&
(1,1,3,1,1)&
(1,2,1,2,1)&
(1,1,1,1,1,1,1). 
\\
\hline
{}[\ld]&\tiny\young(~~~~~~~)
&\tiny\young(~::::,~~~~~,::::~)
&\tiny\young(~~:,~~~,:~~)
&\tiny\young(~~~,:~:,~~~)
&\tiny\young(:~:,:~:,~~~,:~:,:~:)
&\tiny\young(:~:,~~:,:~:,:~~,:~:)
&\tiny\young(~,~,~,~,~,~,~)
\end{array}
\]
\endex
A $\ld$-tableau of type $\B$ is again a bijection
$T: [\ld] \to [d]$. 
We also write $\Sh(T) = \ld$ by a slight abuse of notation.
Such a $\ld$-tableau is called row standard if the entries in each row are increasing, while it is called standard if the entries in each row and each column are increasing.
Denote by $\rStd^\Phi(\ld)$ and $\Std^\Phi(\ld)$ the sets of type $\Phi$ row standard and standard $\ld$-tableaux, respectively. The following assignment is an obvious bijection:
\eq
\rStd^\Phi(\ld) \to \cD^\Phi_\ld,
\quad T \mapsto w_T,
\endeq 
where $w_T$ is given by row-reading of $T$.

Let $\{T\}$ be the equivalence class (called a tabloid) of the $\ld$-tableau $T$ under the (right) action of $W^\A_\ld \subset W^\A_d$.
It is clear that each class contains a unique row standard $\ld$-tableau.

\subsection{Specht module}
We refer to \cite{MR1545531, MR513828, MR644144} for an extensive background on Specht modules. 
For now we assume that $\Phi = \A$.
Define the Specht vector corresponding to $T \in \Sh(\ld)$  by
\eq\label{eq:Specht1}
v_T = \sum_{w \in \textup{col}(T)} (-1)^{\ell^\A(w)} \{T \cdot w\},
\endeq
where $\textup{col}(T)$ is the subset of $\Sigma_d$ consisting of the permutations which reorder the columns of $T$. 
In other words, we have
$\textup{col}(T) = W^\A_{\ld'}$
where $\ld'$ is the transposition of $\ld$, i.e.,
\eq
\ld'_i = \#\{ j \in \ZZ : \ld_j \geq i\}.
\endeq
For $T \in \Std^\Phi(\ld)$, $R \in \rStd^\Phi(\ld)$, and $\ld \in \Lambda^\Phi(d)$, define $c^\A_{R,T}\in\ZZ$ by 
\eq\label{eq:Specht2}
v_T = \sum_{R \in \rStd(d)} c^\A_{R,T} \{R\}.
\endeq

We now generalize this to type B by defining 
\eq\label{eq:SpechtStd}
v_T = \sum_{R \in \rStd^\Phi(d)} c^\Phi_{R,T} \{R\} 
\quad 
\mbox{ for } 
T \in \Sh(\ld)
\mbox{ and } 
\ld \in \Lambda^\Phi(n,d). 
\endeq
\lemma\label{lem:CRT}
If $c^\Phi_{R,T} \neq 0$, then $w_R \leq w_T$ with respect to the Bruhat order.
\endlemma
We denote the (classical) Specht module over $\CC$ corresponding to $\ld \in \Lambda^\Phi(d)$ by 
\eq
S^\ld_\Phi = \bigoplus_{T \in \Std^\Phi(\ld)} \CC v_T,
\endeq
where the (right) action of $W^\Phi_d$ is given by $v^\Phi_T \cdot s_i  = v^\Phi_{T\cdot s_i}$. 
We call $\{v^\Phi_T : T \in \Std^\Phi(\ld)\}$ the tableaux basis for $S^\ld_\Phi$.

Note that for $T \in \Std^\Phi(\ld)$, it is not guaranteed that $T \cdot s_i \in \Std^\Phi(\ld)$, and hence $v^\Phi_{T\cdot s_i}$ is, in general, not a single Specht vector corresponding to a standard $\ld$-tableau, but a linear combination of Specht vectors.
\section{The change-of-basis matrix}
\label{section:change-of-basis}
Now we embed the Specht module $S^J_\Phi$ into the induced trivial module $M^J_\Phi$ by identifying the tabloid $\{R\}$ for $R \in\rStd^\Phi(d)$, with the standard basis element $M^\Phi_{w_R}$ in
$M^J_\Phi$, where $w_R \in \cD^\Phi_J$.

Under the embedding, we denote by $a^\Phi_{x,T}$ (for $T \in \rStd^\Phi(J)$ and $x\in \cD^\Phi_J$) the change-of-basis coefficients from the Specht basis to the Kazhdan--Lusztig basis, i.e.,
\eq
v_T = \sum_{x \in \cD^\Phi_J} a^\Phi_{x,T} M_x.
\endeq
Recall that $\leq$ is the Bruhat order on $\cD^\Phi_J \subset W^\Phi_d$.
\begin{theorem}\label{thm:main1}
Fix total orders $\leq_\cD$ on $\cD^\Phi_J$ and $\leq_r$ on $\rStd^\Phi(J)$ satisfying  
\[
x \leq y \in \cD^\Phi_J \Rightarrow  x \leq_\cD y 
\quad \mbox{ and } 
\quad  
w_R \leq w_T \in \cD^\Phi_J \Rightarrow R \leq_r T, 
\]
respectively. 
The matrix $(a^\Phi_{x,T})_{x,T}$ is upper unitriangular in the sense that
\[
a^\Phi_{w_T, T} =1 \quad \textup{for }T \in \rStd^\Phi(J),
\quad
a^\Phi_{x, T} = 0 \quad \tif x\not\leq_{\cD} w_T.
\]
Moreover, $a^\Phi_{x,T} = 0$ if $x \not\leq w_T$.
\end{theorem}
\proof
Combining Lemmas~\ref{lem:PKL} and \ref{lem:CRT}, we obtain the following:
\eq
v_T 
= \sum_{\substack{R \in \rStd^\Phi(J)\\ w_R \leq w_T}} c^\Phi_{R,T} \sum_{x \leq w_R} m^\Phi_{x, w_R}(1) M_x
= \sum_{x \leq w_T} \sum_{\substack{R \in \rStd^\Phi(J)\\ w_R \leq w_T}} c^\Phi_{R,T}  m^\Phi_{x, w_R}(1) M_x.
\endeq 
It follows from the unitriangularity of $( c^\Phi_{R,T})$ and $(m^\Phi_{w_T, w_R}(1))$ that 
\[
\begin{split}
a^\Phi_{w_T, T} 
=
\sum_{\substack{R \in \rStd^\Phi(J)\\ w_R \leq w_T}} c^\Phi_{R,T}  m^\Phi_{w_T, w_R}(1) 
= 1.
\end{split}\]
Also, if $x\not\leq_{\cD} w_T$ then $x \not\leq w_T$, and so $a^\Phi_{x,T} = 0$.
\endproof
\rmk
For type A, Theorem~\ref{thm:main1} recovers \cite[Theorem~4.1]{Na89} for arbitrary $J$.
If $J = [2n]-\{n\}$, then it also recovers \cite[Theorems~5.5 and 5.7]{RT17}.
\endrmk
\subsection*{Acknowledgment}
The author would like to thank Chun-Ju Lai, Arik Wilbert, and Jieru Zhu for insightful discussions. M.S.I. was partially supported by the Summer Collaborators Program at the School of Mathematics in the Institute for Advanced Study in Princeton, NJ.
\bibliography{litlist} \label{references}
\bibliographystyle{amsalpha}

\end{document}